\documentclass{conm-p-l}

\usepackage{todonotes}
\usepackage{algpseudocode}
\usepackage{pgfplots}
\pgfplotsset{compat=1.16}

\newtheorem{theorem}{Theorem}[section]
\newtheorem{lemma}[theorem]{Lemma}
\newtheorem{corollary}[theorem]{Corollary}

\theoremstyle{definition}
\newtheorem{definition}[theorem]{Definition}
\newtheorem{example}[theorem]{Example}

\theoremstyle{remark}
\newtheorem{remark}[theorem]{Remark}

\numberwithin{equation}{section}

\newcommand{\gdag}{g^\dagger}
\newcommand{\gobs}{g^{\mathrm{obs}}} 
\newcommand{\udag}{u^\dagger}
\newcommand{\X}{\mathcal{X}}
\newcommand{\Y}{\mathcal{Y}}
\newcommand{\N}{\mathbb N}
\newcommand{\R}{\mathbb R}

\newcommand{\Yd}{\mathcal Y_{d}}

\renewcommand{\L}{\mathbf L}
\renewcommand{\H}{\mathbf H}
\newcommand{\app}{\mathrm{app}}

\newcommand{\W}{\mathbf W}
\newcommand{\manifold}{\mathbb M}
\newcommand{\err}{\mathrm{err}}
\newcommand{\norm}[2]{\left\Vert #1 \right\Vert_{#2}}
\newcommand{\E}{\mathbb{E}}
\DeclareMathOperator*{\argmin}{argmin}
	\newcommand{\res}{\mathrm{res}}

\newcommand{\abs}[1]{\lvert#1\rvert}

\newcommand{\new}[1]{#1}

\newlength\fheight \newlength\fwidth % Power plots

\begin{document}

\title[Optimal rates for $\mathbf L^1$ data fitting]{$\mathbf L^1$ data fitting for Inverse Problems yields optimal rates of convergence in case of discretized white Gaussian noise}

%    Information for first author
\author{Kristina Bätz}
%    Address of record for the research reported here
\address{Institute of Mathematics, Julius-Maximilians-Universität Würzburg, 97084 Würzburg}
\email{kristina.baetz@uni-wuerzburg.de}
%    \thanks will become a 1st page footnote.
\thanks{This manuscript reports on the PhD thesis of the first author.}

%    Information for second author
\author{Frank Werner}
\address{Institute of Mathematics, Julius-Maximilians-Universität Würzburg, 97084 Würzburg}
\email{frank.werner@uni-wuerzburg.de}

%    General info
\date{\today}


\begin{abstract}
It is well-known in practice that $\L^1$ data fitting leads to improved robustness compared to standard $\L^2$ data fitting. However, it is unclear whether resulting algorithms will perform as well in case of regular data without outliers. In this paper, we therefore analyze generalized Tikhonov regularization with $\L^1$ data fidelity for Inverse Problems $F(u) = g$ in a general setting, including general measurement errors and errors in the forward operator. The derived results are then applied to the situation of discretized Gaussian white noise, and we show that the resulting error bounds allow for order-optimal rates of convergence. These findings are also investigated in numerical simulations. 
\end{abstract}

\maketitle

\section{Introduction}

In this paper we study Inverse Problems
\begin{equation}\label{eq:ip}
 F(u) = g
\end{equation}
with an operator $F : D(F) \subset \X\to \L^1 \left(\manifold\right)$ mapping between a Banach space $\X$ and $\L^1 \left(\manifold\right)$, where $\manifold \subset \R^d$ is some smooth bounded domain. The exact solution of \eqref{eq:ip} will be denoted by $\udag \in \X$, and the corresponding exact datum by $\gdag := F(\udag) \in \L^1 \left(\manifold\right)$. As typical in applications, we will assume that $\gdag$ is not (exactly) available, but rather some approximation $\gobs \in \L^1 \left(\manifold_h\right)$ of it, where $\left(\manifold_h,\Sigma,\nu\right)$ is a (possibly discrete) finite measure space. Note that this allows $\gobs$ to live on a different domain than $\gdag$, e.g. for discretization reasons. This setup also allows for finite-dimensional data in the sense that $\gobs$ can be interpreted as a vector in $\R^n$, but also for different types of infinite-dimensional data. For the moment we do not specify a model for $\gobs$. Additionally, we will assume that $F$ cannot be evaluated or is not known exactly, but only some approximation
\[
F_h : D(F) \to \L^1 \left(\manifold_h\right)
\]
is at hand. Here and above, $h > 0$ can be seen as a discretization or approximation parameter. One example for the approximations $F_h$ is given by sampling or taking local averages, which can formally be written as $F_h = P_h F$ for some operators $P_h : \L^1 \left(\manifold\right)\to \L^1 \left(\manifold_h\right)$. However, our approximation model also includes the case of a numerical approximation of the forward operator $F$ e.g. by finite element schemes in case that $F$ is described by means of partial differential equations.

If the available data $\gobs$ might be corrupted by outliers, a common approach to reconstruct $\udag$ from $\gobs$, knowing only $F_h$, is given by (generalized) Tikhonov regularization with $\L^1$ data fidelity, this is
\begin{equation}\label{eq:tik}
	\hat u_\alpha \in \argmin_{u \in D(F)} \left[\norm{F_h\left(u\right)- \gobs}{\L^1 \left(\manifold_h\right)}+ \alpha \mathcal R \left(u\right)\right],
\end{equation}
where $\mathcal R : \X \to \left(-\infty, \infty\right]$ is a regularizer or penalty functional, and $\alpha > 0$ is a regularization parameter. Typically, $\mathcal R$ is chosen to be convex and to obey some continuity properties w.r.t. a sufficiently weak topology on $\X$. 

The rationale behind the choice of the $\L^1$ data fidelity term $\norm{\cdot}{\L^1 \left(\manifold_h\right)}$ in \eqref{eq:tik} is, that compared to the standard least-squares choice $\norm{\cdot}{\L^2\left(\manifold_h\right)}^2$, it is somewhat more robust w.r.t. outliers in the data. If e.g. $y_1,\ldots,y_n \in \R$ are observations, then the minimizer $\bar y := \argmin_{y \in \R} \sum_{i=1}^n \left(y-y_i\right)^2$ is the empirical mean $\bar y = n^{-1}\sum_{i=1}^n y_i$. If one of the observations is far \new{off} the others, then this influences $\bar y$ with a weight $n^{-1}$. A more robust estimator is the median, which is given by $\tilde y := \argmin_{y \in \R} \sum_{i=1}^n \left|y-y_i\right|$, and is well-known to be independent of the specific value of the smallest or largest $y_i$, since it can practically be computed as the $1/2$-quantile of $y_1,\ldots,y_n$. As a consequence, a better choice as data fidelity term in case of possible outliers in the data $\gobs$ is $\norm{\cdot}{\L^1\left(\manifold_h\right)}$ as in \eqref{eq:tik}. This yields arguably a robust method, see e.g. \cite{hw14,kwh16,cj12} and the references therein. 

Generalized Tikhonov regularization, i.e. methods of the form \eqref{eq:tik}, \new{has} been studied intensively over the last decades, see e.g. \cite{e93,bo04,r05,hkps07,bkmss08,g10b,bb11,wh12,f18b,kr19,bbks20,hp20,ph25} to name a few.
Investigated questions include existence of minimizers, possible uniqueness, stability w.r.t. $\gobs$, and finally also convergence $\hat u_\alpha \to \udag$ if $\gobs \to \gdag$, $F_h \to F$ and $\alpha \searrow 0$ in a suitable sense.

In this paper, we will mostly be interested in rates of convergence, \new{that} is\new{,} quantitative estimates of $\ell \left(\hat u_\alpha, \udag\right)$ for some loss function $\ell : \X \times \X\to \left[0,\infty\right]$ in terms of $\alpha$, the noise contained in $\gobs$, the approximation parameter $h$ of $F_h$, and smoothness properties of $\udag$. Since our final model for $\gobs$ will be random, we are also interested in the expectation or concentration of $\ell\left(\hat u_\alpha,\udag\right)$. Once such error bounds for $\ell\left(\hat u_\alpha,\udag\right)$ are derived, a natural question is the optimality of those, i.e. whether methods different from \eqref{eq:tik} such as using $\L^2$ data fitting instead or different proof techniques could lead to improved estimates. 

For many specific classes of operators \new{$F : D(F) \subseteq \L^2 \left(\Omega\right) \to \L^2 \left(\manifold\right)$} and noise models, the optimal rates of convergence are known, see e.g. \cite{bhmr07,c08,rsh16,wsh20}. Especially, if the noise is white Gaussian with variance $\sigma^2> 0$, and if the operator is a linear isomorphism between $H^t \left(\Omega\right)$ and $H^{t+a}\left(\manifold\right)$ for all $t \in \R$ with some fixed $a> 0$, where $\Omega \subset \R^m$ is another smooth bounded domain, then under a source condition $\udag \in H^s \left(\Omega\right)$ with $0 < s < a$, the optimal rate of convergence is
\begin{equation}\label{eq:rate_opt}
\left(\mathbb E \left[ \norm{\hat u_\alpha - \udag}{\L^2 \left(\Omega\right)}^2\right]\right)^{\frac12} = \mathcal O \left(\sigma^{\frac{2s}{2a+2s+d}}\right),
\end{equation}
see e.g. \cite[Thm. 4.10]{wsh20}. \new{Note that any such operator can, by boundedness of $\manifold$ and hence $\L^\infty \left(\manifold\right) \hookrightarrow \L^1 \left(\manifold\right)$, also be considered as an operator mapping into $\L^1 \left(\manifold\right)$ as soon as $a > d/2$ for $\X = \L^2 \left(\Omega\right)$. This is the case which will be considered mostly below.}

In this paper, we aim to provide a new analysis for generalized Tikhonov regularization with $\L^1$ data fidelity \eqref{eq:tik}. If $\gobs$ is given as $\gdag$ plus impulsive noise (see the definition in Section 3 below), which is a continuous model for noise containing arbitrarily heavy outliers, this has been analyzed in \cite{hw14,kwh16,hw22}. We build upon those results and aim to answer what happens if \eqref{eq:tik} is used in a discretized observation scheme. We will show that an extension of the previously mentioned analysis can be used to derive rates of convergence, which are known to be optimal in the Gaussian case. We also discuss efficient numerical methods for the minimization of \eqref{eq:tik} and investigate their performance in simulations. As a consequence, we recommend to use $\L^1$ data fitting as the standard method, since it has nearly no drawbacks compared to $\L^2$ data fitting, but provides a \new{significantly} more robust algorithm.

The remainder of this paper is organized as follows. In Section 2 we recall central results on generalized Tikhonov regularization \eqref{eq:tik}, especially the tools required to derive error bounds as presented in \cite{hw22}. The \new{successive} Section 3 is devoted to impulsive noise and provides the specific model introduced in \cite{hw14,kwh16} and the resulting convergence rates in this case. In Section 4 we discuss the more realistic setting of discrete noise and derive error bounds for $\ell\left(\hat u_\alpha,\udag\right)$, especially in case of discretized white Gaussian noise. The resulting rates of convergence turn out to equal \eqref{eq:rate_opt}, i.e. are known to be optimal in this case. Section 5 presents numerical simulations, which support the previous estimates. In Section 6 we briefly discuss our results and conclude.

\section{Error bounds for generalized Tikhonov regularization}

Let us briefly recall the setup from \cite{hw22}, where error bounds for generalized Tikhonov regularization \eqref{eq:tik} are provided. 

\subsection{Effective noise level}

To measure the influence of the errors $F_h \approx F$ and $\gobs \approx \gdag$, the so-called \textbf{effective noise level functional} $\err : D(F) \to (-\infty,\infty)$ is introduced by
\[
\err(u) := \norm{F(u)- \gdag}{\L^1\left(\manifold\right)} - \left(\norm{F_h(u)-\gobs}{\L^1 \left(\manifold_h\right)} - \norm{F_h \left(\udag\right)-\gobs}{\L^1 \left(\manifold_h\right)}\right).
\]
The rationale behind is that we would like to use $u \mapsto \norm{F(u)- \gdag}{\L^1\left(\manifold\right)}$ as a data fidelity, i.e. to compute
\begin{equation}\label{eq:tik_exact}
		u_\alpha \in \argmin_{u \in D(F)} \left[\norm{F\left(u\right)- \gdag}{\L^1 \left(\manifold\right)}+ \alpha \mathcal R \left(u\right)\right],
\end{equation}
corresponding to the exact operator and noise\new{-}free data, but this is not possible due to lack of knowledge. The definition of $\err$ implies that
\[
\norm{F \left(\hat u_\alpha\right)-\gdag}{\L^1 \left(\manifold\right)}+ \alpha \mathcal R \left(\hat u_\alpha\right) \leq \norm{F\left(\udag\right)-\gdag}{\L^1 \left(\manifold\right)} + \alpha\mathcal R \left(\udag\right) + \err \left(\hat u_\alpha\right)
\]
for any minimizer $\hat u_\alpha$ as in \eqref{eq:tik}, i.e. $\hat u_\alpha$ somewhat minimizes the ideal Tikhonov functional \eqref{eq:tik_exact} up to $\err \left(\hat u_\alpha\right)$. This means that all error contributions from $F_h \approx F$ and $\gobs \approx \gdag$ are factored out into $\err \left(\hat u_\alpha\right)$. If this term can be bounded, this allows for deriving error bounds.

\subsection{Variational source condition}

The error contribution of (ideal) generalized Tikhonov regularization \eqref{eq:tik_exact} can be controlled as follows. Suppose we want to investigate the error w.r.t. some \textbf{loss} function $\ell : \X \times \X \to \left[0,\infty\right]$. 
\new{The m}ost common examples of losses are norm powers $\ell \left(u,\udag\right)= \norm{u-\udag}{\X}^r$ or Bregman divergences
\[
\ell \left(u,\udag\right) = \mathcal D_{\mathcal R}^{u^*} \left(u;\udag\right) = \mathcal R(u) - \mathcal R\left(\udag\right) - \left\langle u^*, u-\udag\right\rangle
\]
where $u^* \in \partial \mathcal R \left(\udag\right)$ \new{denotes an element of the subdifferential $\partial \mathcal R \left(\udag\right)$}. Now let us assume that $\udag \in D(F)$ satisfies a \textbf{variational source condition}
\begin{equation}\label{eq:vsc}
	\ell\left(u,\udag\right) \leq \mathcal R(u) - \mathcal R \left(\udag\right) + \varphi \left(\norm{F(u)-\gdag}{\L^1 \left(\manifold\right)}\right) \qquad\text{for all}\qquad u \in D(F)
\end{equation}
with an \textbf{index function} $\varphi:[0,\infty) \to [0,\infty)$, i.e. $\varphi$ is continuous, increasing, and obeys $\varphi(0) = 0$. 

If we now compare the minimizer $u_\alpha$ from \eqref{eq:tik_exact} with $\udag$, we obtain using \eqref{eq:vsc} that
\[
\ell\left(u_\alpha;\udag\right) \leq -\frac1\alpha\norm{F \left(u_\alpha\right)-\gdag}{\L^1 \left(\manifold\right)} + \varphi \left(\norm{F \left(u_\alpha\right)-\gdag}{\L^1 \left(\manifold\right)}\right) \leq \varphi_\app \left(\alpha\right)
\]
with the convex conjugate $\varphi_\app \left(\alpha\right) := \sup_{\tau \geq 0} \left[-\frac1\alpha \tau + \varphi\left(\tau\right)\right] = \left(-\varphi\right)^* \left(-\frac1\alpha\right)$, see \cite{hw22}. This shows that the approximation error introduced by (generalized) Tikhonov regularization can be bounded in terms of $\varphi_{\app}$, a function which has natural properties (see \cite[Rem. 2]{wh20}, \cite[Rem. 3.29]{w12} and \cite{hw22}):
\begin{enumerate}
	\item[(a)] $\varphi_\app \left(\alpha\right) \geq 0$ for all $\alpha > 0$;
	\item[(b)] $\varphi_\app$ is monotonically increasing; and
	\item[(c)] if $\varphi^2$ is concave, then $\lim_{\alpha \searrow 0}\varphi_{\app}(\alpha) = 0$. %; and
%	\item[(d)] if $\varphi^(1+\epsilon)$ is concave for some $\epsilon > 0$, then
%	\[
%	\varphi_{\app}\left(C\alpha\right) \leq C^{\frac1\epsilon} \varphi_{\app} \left(\alpha\right) \qquad\text{for all} \qquad C \geq 1, \alpha > 0.
%	\]
\end{enumerate}

Variational source conditions \eqref{eq:vsc} have been introduced in \cite{hkps07}, and have since then developed as a central tool for the analysis of variational regularization, see e.g. \cite{f18a,f18b}. They have been verified using different techniques for specific examples of Inverse Problems, see e.g. \cite{hw15,hw17,wh17,w19}. 

\begin{example}\label{ex:vsc}
	As an illustrative example, assume that $\Omega \subset \R^m$ and $\manifold \subset \R^d$ are smooth bounded domains, and $F : D(F) \subset \L^2 \left(\Omega\right) \to \L^2 \left(\manifold\right)$ is \textbf{at most $a$-times smoothing} in the sense that
	\begin{equation}\label{eq:a_smoothing_most}
	\norm{u-v}{\H^{-a} \left(\Omega\right)} \leq C_F \norm{F(u) - F(v)}{\L^2 \left(\manifold\right)} \qquad\text{for all}\qquad u,v \in D(F),
	\end{equation}
	where $a > 0$ is a parameter. Then the classical smoothness condition $\udag \in \H^s_0 \left(\Omega\right)$ for some $0 < s < a$ with $\norm{\udag}{\H^s} \leq \rho$ implies the variational source condition 
	\begin{equation}\label{eq:vsc_L2}
	\frac14 \norm{u-\udag}{\L^2\left(\manifold\right)}^2 \leq \frac12 \norm{u}{\L^2\left(\Omega\right)}^2 - \frac12 \norm{\udag}{\L^2 \left(\Omega\right)} +  \varphi \left(\norm{F(u)-\gdag}{\L^2\left(\manifold\right)}^2\right)
	\end{equation}
	with $\varphi(t) = C \rho^{\frac{2a}{a+s}} t^{\frac{s}{a+s}}, t > 0$, for some constant $C > 0$ independent of $\rho$ and $t$, see Example 2.2 in \cite{hw22}. In this case, it follows via \new{convex conjugation properties} that
\[
\varphi_\app \left(\alpha\right) = \tilde C \rho^2 \alpha^{\frac{s}{a}}
\]
with some different constant $\tilde C > 0$ again independent of $\rho$ and $\alpha$. 
\end{example}

For a possibility to derive a variational source condition \eqref{eq:vsc} with $\L^1$ data fidelity from \eqref{eq:vsc_L2}, see \cite[Prop. 2.2]{kwh16}.

\subsection{Bounding the effective noise level}

To further bound the effective noise level $\err\left(\hat u_\alpha\right)$, two further concepts have been introduced in \cite{hw22}: \textbf{admissibility} and \textbf{smoothingness}. 

Admissibility assumes an abstract upper bound
\begin{equation}\label{eq:admissibility}
\err \left(u\right) \leq \eta_0 + \sum_{m=1}^M \eta_m \norm{F(u) - \gdag}{\Y_{\mathrm{d}}^m}, \qquad u \in D(F)
\end{equation}
for the effective noise level functional, \new{reflecting} the fact that, especially for random and / or operator noise, this functional will not be uniformly bounded for all $u \in D(F)$. In \eqref{eq:admissibility}, $\Y_{\mathrm{d}}^m \subset \L^1 \left(\manifold\right)$ for $m=1,\ldots,M$ are Banach spaces with $R(F) \subset \Y_{\mathrm{d}}^m$, and $\eta_0,\ldots, \eta_M \geq 0$ are random variables (or constants in the deterministic case). Any specific model for $\gobs$ and $F_h$ in \eqref{eq:ip} satisfying \eqref{eq:admissibility} is called $\left(\Y_{\mathrm{d}}^1, \ldots,\Y_{\mathrm{d}}^M\right)$-admissible.

The additional residuals $ \norm{F(u) - \gdag}{\Y_{\mathrm{d}}^m}$ in \eqref{eq:admissibility} are then treated by interpolation and smoothing properties of $F$. The operator $F$ is called $\left(\Y_{\mathrm{d}}, \gamma, \ell\right)$-smoothing, if
\begin{equation}\label{eq:smoothingness}
\norm{F(u) - \gdag}{\Y_{\mathrm{d}}} \leq \gamma \left(\delta\right) \ell \left(u,\udag\right)^{\frac1r} + \frac1\delta \norm{F(u) - \gdag}{\L^1 \left(\manifold\right)}, \qquad u \in D(F),
\end{equation}
for all $\delta \in [0,\delta_0]$, where $\gamma : [0,\delta_0] \to [0,\infty)$ is a function and $\delta_0 > 0$, $r \geq 1$ are constants. The parameter $r$ is only to take care of potential exponents in $\ell$ e.g. in case of norm losses. The bound \eqref{eq:smoothingness}, if assumed for each space $\Y_{\mathrm{d}}^m$, $m=1,\ldots,M$, allows \new{us} to further estimate each of the additional terms in \eqref{eq:admissibility} separately.

It has been argued in \cite{hw22}, that \eqref{eq:smoothingness} is in fact a combination of smoothing properties of $F$ with interpolation properties of $\L^1 \left(\manifold\right)$ and $\Y_{\mathrm{d}}$, cf. Remark 2.7 there. Specifically for $\Y_{\mathrm{d}} = \L^\infty \left(\manifold\right)$, this can be shown explicitly:
\begin{lemma}\label{lem:smoothingness}
Suppose $\manifold \subset \R^d$ is a bounded Lipschitz domain and $F : D(F) \subset X \to \W^{k,p} \left(\manifold\right)$ is Lipschitz continuous in the sense that
\begin{equation}\label{eq:lip}
	\norm{F(u) - F\left(\udag\right)}{\W^{k,p} \left(\manifold\right)} \leq L \ell \left(u,\udag\right)^{\frac1r},
\end{equation}
\new{where $k > \frac{d}{p}$ and $\W^{k,p} \left(\manifold\right)$ denotes the standard $\L^p$-based Sobolev space.}
Then $F$ is $\left(\L^\infty\left(\manifold\right), \gamma,\ell\right)$-smoothing with $\gamma\left(\delta\right) = C \delta ^{\frac{k}{d} - \frac1p}$ for some constant $C > 0$ independent of $\delta$.
\end{lemma}
\begin{proof}
This estimate is derived in \cite[Sec. 6]{hw22}, and is based on an application of Ehrling's lemma (cf. \cite{hw14} for details).
\end{proof}

\begin{remark}
If $F$ is \textbf{at least $a$ smoothing} in the sense that
\begin{equation}\label{eq:a_smoothing_least}
	\norm{F(u) - F\left(v\right)}{\H^a \left(\manifold\right)} \leq C_F \norm{u-v}{\L^2 \left(\Omega\right)}\qquad \text{for all}\qquad u,v \in D(F),
\end{equation}
where $a > 0$ is again a parameter, then\new{,} for $a \in \N$ and $\ell = \norm{\cdot - \cdot}{\L^2 \left(\Omega\right)}$\new{,} this implies \eqref{eq:lip} with $k = a$ and $p =2$. In this sense, \eqref{eq:lip} can be seen as an at least smoothing assumption compared to the at most smoothing assumption \eqref{eq:a_smoothing_most}. 
\end{remark}

\subsection{Error bounds}

All together, we obtain the following error bounds for generalized Tikhonov regularization \eqref{eq:tik} with $\L^1$ data fidelity (see \cite[Thm. 3.4]{hw22}):
\begin{theorem}\label{thm:bounds}
Suppose the data model for $\gobs$ is $\left(\Y_{\mathrm{d}}^1,\ldots,\Y_{\mathrm{d}}^M\right)$-admissible with $\eta_0,\ldots,\eta_M \geq 0$ as in \eqref{eq:admissibility} and that the variational source condition \eqref{eq:vsc} is satisfied. \new{Suppose furthermore that $r > 1$.} For those $m \in \{1,\ldots,M\}$ with $\eta_m > 0$, we additionally assume that $F$ is $\left(\Y_{\mathrm{d}}^m, \gamma_m, \ell\right)$-smoothing as in \eqref{eq:smoothingness}. If $4 \eta_m < \delta_0^m$ for all such $m$, then each global minimizer $\hat u_\alpha$ of \eqref{eq:tik} (if \new{it exists}) satisfies
\begin{equation}
\ell\left(\hat u_\alpha,\udag\right) \leq C \left[\frac{\eta_0}{\alpha} + \sum_{\genfrac{}{}{0pt}{}{m=1,\ldots,M}{\eta_m > 0}} \left(\frac{\eta_m}{\alpha}\right)^{r'} \left(\gamma_m \left(\eta_m\right)\right)^{r'}+ \varphi_{\app} \left(2\alpha\right)\right]
\end{equation}
with the conjugate exponent $r' := \frac{r}{r-1}$ of $r$, and a constant $C > 0$ independent of $\eta_0, \ldots, \eta_M$ and $\alpha$.
\end{theorem}

\section{Impulsive noise and resulting error bounds}

In this section, let us now specify the results from Section 2 on generalized Tikhonov regularization with $\L^1$ data fitting \eqref{eq:tik} to the situation of impulsive noise.

\subsection{Impulsive noise}

In general, a noise vector or function $\xi: \manifold_h \to \R$ is called \textbf{impulsive} if $\abs{\xi}$ is large on a small part of $\manifold_h$, and small or even zero elsewhere. This type of noise occurs e.g. in measurements with malfunctioning receivers, digital image acquisition with faulty memory locations, or in powerline communication systems, see also \cite{chn04,cj12}. 

In \cite{hw14}, a continuous model for impulsive noise has been introduced as follows:
\begin{definition}[Impulsive noise]
	Let $\left(\manifold_h, \Sigma,\nu\right)$ be a (possibly discrete) measure space on which the observations \new{are defined}. A noise function $\xi \in \L^1\left(\manifold_h\right)$ is called \textbf{impulsive} with parameters $\varepsilon, \eta \geq 0$, if there exists $\manifold_{\mathrm{c}} \in \Sigma$ such that
	\[
	\norm{\xi}{\L^1 \left(\manifold_h \setminus \manifold_{\mathrm{c}}\right)} \leq \varepsilon, \qquad \nu \left(\manifold_{\mathrm{c}}\right) \leq \eta.
	\]
\end{definition}
In the above setting, the subset $\manifold_{\mathrm{c}} \subset \manifold_h$ is considered as the corrupted part of the measurements, whereas $\manifold_{\mathrm{i}} := \manifold_h \setminus \manifold_{\mathrm{c}}$ is considered as the intact part where the noise is small.

In the following, we will consider the noise model
\begin{equation}\label{eq:noise_model}
\gobs := F_h \left(\udag\right) + \xi
\end{equation}
with an impulsive noise function $\xi = \xi \left(\varepsilon,\eta\right) \in \L^1 \left(\manifold_h\right)$. Note that we do not require any (direct) relation between $\gobs$ and $\gdag = F \left(\udag\right)$ here, but only between $\gobs$ and the natural discretization $F_h \left(\udag\right)$. The misfit between $F_h \left(\udag\right)$ and $\gdag$ will be bounded \new{in the subsection} below by assuming properties of $F_h$ as an approximation of $F$.

\subsection{Resulting error bounds}

Let us now specify the error bounds derived in Theorem \ref{thm:bounds} for the case of impulsive noise as discussed above along the lines of \cite{hw22}. To do so, we will assume that the operator approximation $F_h : D(F) \subset \X \to \L^1 \left(\manifold_h\right)$ satisfies the error bounds
\begin{subequations}\label{eqs:discretization}
	\begin{align} 
		\norm{F\left(u\right)- \gdag}{\L^1\left(\manifold\right)}- \norm{F_h\left(u\right)-F_h\left(\udag\right)}{\L^1\left(\manifold_h\right)} &\leq C_\omega h^\omega \norm{F\left(u\right)-\gdag}{\Y_{\mathrm{d}}},\label{eq:discretization1} \\
		\norm{F_h\left(u\right)-F_h\left(\udag\right)}{\L^\infty\left(\manifold_h\right)} &\leq C_\omega \norm{F\left(u\right)-\gdag}{\L^\infty\left(\manifold\right)}\label{eq:discretization2}
	\end{align}
\end{subequations}
for all $u\in D(F)$ with constants $C_\omega, \omega\geq 0$ and a norm $\norm{\cdot}{\Y_{\mathrm{d}}}$.

Assumptions of the form \eqref{eqs:discretization} can be verified for different approximations $F_h$. As discussed in Section 4 below, it is almost trivial for situations where $F_h$ corresponds to a point-sampled version of $F$. However, \eqref{eqs:discretization} can also be verified in more general settings where other sampling schemes are used or $F_h$ is directly defined by means of finite element methods, see e.g. 
\cite{o94,dvp88,bs94}.

We will use \cite[Eq. (3.1)]{hw14} for our analysis below. There, it has been shown that the estimate 
\begin{equation}\label{eq:aux}
	\begin{split}
\norm{g- F_h \left(\udag\right)}{\L^1\left(\manifold_h\right)} - &\left(\norm{g-\gobs}{\L^1 \left(\manifold_h\right)} - \norm{\xi}{\L^1 \left(\manifold_h\right)}\right) \\&\leq 2 \norm{\xi}{\L^1 \left(\manifold_h \setminus \manifold_{\mathrm{c}}\right)} + 2 \nu \left(\manifold_{\mathrm{c}}\right)\norm{g - F_h \left(u\right)}{\L^\infty \left(\manifold_h\right)}\\
& \leq 2 \varepsilon + 2 \eta \norm{g - F_h \left(u\right)}{\L^\infty \left(\manifold_h\right)}
\end{split}
\end{equation}
holds true for any $g \in \L^1 \left(\manifold_h\right)$ and any impulsive noise function $\xi = \xi \left(\varepsilon,\eta\right) \in \L^1 \left(\manifold_h\right)$. Using this, we obtain the following result:
\begin{theorem}\label{thm:admissibility}
	Suppose the noise model \eqref{eq:noise_model} with an impulsive noise $\xi = \xi \left(\varepsilon,\eta\right) \in \L^1 \left(\manifold_h\right)$, and let $F_h$ satisfy \eqref{eqs:discretization}. Then the overall data model is $\left(\L^\infty\left(\manifold\right), \Y_{\mathrm{d}}\right)$-admissible \new{(i.e. $\Y_{\mathrm{d}}^1 = \L^\infty\left(\manifold\right)$ and $\Y_{\mathrm{d}}^2 = \Y_{\mathrm{d}}$)} with $\eta_0 = 2 \epsilon$, $\eta_1 = 2 C_\omega\eta $ and $\eta_2=C_\omega h^\omega$.
\end{theorem}
\begin{proof}
For all $u \in D(F)$ we have
\begin{align*}
	\err\left(u\right) 
	= &\norm{F(u)-\gdag}{\L^1 \left(\manifold\right)}- \left(\norm{F_h(u) - \gobs}{\L^1 \left(\manifold_h\right)}- \norm{F_h\left(\udag\right)-\gobs}{\L^1 \left(\manifold_h\right)}\right) \\
	= &\norm{F(u)-\gdag}{\L^1 \left(\manifold\right)}- \left(\norm{F_h(u) - \gobs}{\L^1 \left(\manifold_h\right)}- \norm{\xi}{\L^1 \left(\manifold_h\right)}\right) \\
	\leq & \norm{F_h\left(u\right)-F_h\left(\udag\right)}{\L^1\left(\manifold_h\right)} + C_\omega h^\omega \norm{F\left(u\right)-\gdag}{\Y_{\mathrm{d}}}\\
	& - \left(\norm{F_h(u) - \gobs}{\L^1 \left(\manifold_h\right)}- \norm{\xi}{\L^1 \left(\manifold_h\right)}\right) \\
	\leq &2 \varepsilon + 2\eta \norm{F_h(u)-F_h\left(\udag\right)}{\L^\infty\left(\manifold_h\right)} + C_\omega h^\omega \norm{F\left(u\right)-\gdag}{\Y_{\mathrm{d}}}\\
	\leq &2 \varepsilon + 2\eta C_\omega \norm{F(u)-F\left(\udag\right)}{\L^\infty\left(\manifold\right)} + C_\omega h^\omega \norm{F\left(u\right)-\gdag}{\Y_{\mathrm{d}}},	
	\end{align*}
	where we used \eqref{eq:discretization1} first, then \eqref{eq:aux} with $g = F_h(u)$, and finally \eqref{eq:discretization2}. This shows the claimed admissibility.
\end{proof}

With this estimate, we are now in \new{a} position to specify the result from Theorem \ref{thm:bounds}, cf. also \cite[Cor. 6.1]{hw22}
\begin{corollary}\label{cor:bounds}
	Suppose the noise model \eqref{eq:noise_model} with an impulsive noise $\xi = \xi \left(\varepsilon,\eta\right) \in \L^1 \left(\manifold_h\right)$, and let $F_h$ satisfy \eqref{eqs:discretization}. Let furthermore the variational source condition \eqref{eq:vsc} hold true, $F: D(F) \to \W^{k,p}\left(\manifold\right)$ be Lipschitz continuous with $k > \frac{d}{p}$ as in \eqref{eq:lip} and suppose $F$ is $\left(\Yd, \gamma,\ell\right)$-smoothing. \new{Assume $r > 1$.} Then there exists a constant $C$ independent of $\alpha, \varepsilon, \eta$ and $h$ such that any minimizer of \eqref{eq:tik} (if \new{it exists}) satisfies the error bound
	\begin{align}\label{eq:err_bound}
		\ell\left(\hat u_\alpha,\udag\right) \leq C \left[\frac{\varepsilon}{\alpha} + \left(\frac{\eta^{\vartheta+1}}{\alpha}\right)^{r'}+ \left(\frac{h^\omega}{\alpha}\gamma \left( h^\omega\right)\right)^{r'}+ \varphi_{\app}(2\alpha)\right],
	\end{align}
where again $r' = \frac{r}{r-1}$ and $\vartheta: =\frac{k}{d}-\frac{1}{p} > 0$.
\end{corollary}

\begin{remark}
	Note that the above result can be generalized to exponentially smoothing operators, i.e. to the case that $F: D(F) \to \W^{k,p}\left(\manifold\right)$ is Lipschitz for any choice of $k \in \N$, see \cite{kwh16}. There, similar bounds with $\eta^{\vartheta+1}$ replaced by $\tilde \gamma \left(\eta\right)$ have been derived, where the function $\tilde \gamma$ tends faster to zero than any polynomial. For the backwards heat equation, one e.g. obtains $\tilde \gamma \left(\delta\right) \sim \exp\left(-\frac{1}{\delta^2}\right)$, and for three-dimensional satellite gradiometry $\tilde \gamma \left(\delta\right) \sim \delta^{-\frac52} R^{-\sqrt{\frac{4\pi}{\delta}-4}}$, where $R > 1$ is the radius of the measurement shell. However, for simplicity, we restrict ourselves to the finitely smoothing case as discussed above in this work.
\end{remark}
\subsection{Results for finitely smoothing operators under classical smoothness conditions}

Let us now derive error bounds in the setting of Example \ref{ex:vsc} for generalized Tikhonov regularization with $\L^1$ data fidelity term. Note that Corollary \ref{cor:bounds} is not directly applicable, since Example \ref{ex:vsc} only implies the variational source condition \eqref{eq:vsc_L2} with $\L^2$ data fidelity term. We therefore exploit ideas from \cite[Prop. 2.2]{kwh16} and mimic the proofs from \cite{hw22} to obtain the following:
\begin{theorem}\label{thm:example}
	Let $\Omega \subset \R^m$ and $\manifold \subset \R^d$ be smooth bounded domains, and choose
	\[
	\mathcal R \left(u\right) := \frac12 \norm{u}{\L^2 \left(\Omega\right)}^2,\qquad \ell \left(u,\udag\right):= \frac14 \norm{u-\udag}{\L^2 \left(\Omega\right)}^2.
	\]
	Furthermore let $F : D(F) \subset \L^2 \left(\Omega\right) \to \L^2 \left(\manifold\right)$ be the forward operator. Suppose the following conditions hold true:
	\begin{enumerate}
		\item[(A1)] The noise model is \eqref{eq:noise_model} with impulsive noise $\xi = \xi \left(\varepsilon,\eta\right) \in \L^1 \left(\manifold_h\right)$.
		\item[(A2)] $F$ is at most $a$-times smoothing in the sense that \eqref{eq:a_smoothing_most} holds true for all $u, v \in D(F)$ with a constant $C_F > 0$ and parameter $a > 0$.
		\item[(A3)] $\udag \in \H_0^s \left(\Omega\right)$ for some $0 < s < a$ with $\norm{\udag}{\H^s\left(\Omega\right)} \leq \rho$.
		\item[(A4)] $F : D(F) \subset \L^2 \left(\Omega\right) \to \W^{k,p}\new{\left(\manifold\right)}$ is Lipschitz continuous w.r.t. $\ell$ where $k > \frac{d}{p}$, th\new{at} is\new{,}
		\begin{equation}\label{eq:lip_norm}
			\norm{F(u) - F\left(\udag\right)}{\W^{k,p} \left(\manifold\right)} \leq L \norm{u-\udag}{\L^2 \left(\Omega\right)}
		\end{equation}
		for all $u \in D(F)$ with a constant $L > 0$.
		\item[(A5)] The approximation $F_h$ of $F$ satisfies \eqref{eqs:discretization} with some norm $\norm{\cdot}{\Y_{\mathrm{d}}}$ with $\omega > 0$ and a constant $C_\omega \geq 0$.
		\item[(A6)] $F$ is $\left(\Y_{\mathrm{d}}, \gamma, \norm{\cdot - \cdot}{\L^2 \left(\Omega\right)}^2\right)$-smoothing, where $r = 2$.
	\end{enumerate}
	Then, provided $\alpha, \varepsilon, \eta$ and $h$ are sufficiently small, any minimizer $\hat u_\alpha$ of \eqref{eq:tik} (if \new{it exists}) obeys the error bound
	\begin{equation}\label{eq:bound_example}
		\norm{\hat u_\alpha -\udag}{\L^2 \left(\Omega\right)}^2 \leq C \left[\frac{\varepsilon}{\alpha} + \frac{\eta^{2\vartheta}}{\alpha^2}+\frac{h^{2\omega}\left(\gamma \left( 3 h^\omega\right)\right)^2}{\alpha^2}+ \alpha^{\frac{2s(2\vartheta+1)}{(2\vartheta+1)(a-s)+a}}\right]
	\end{equation}
	where $\vartheta = \frac{k}{d} - \frac1p > 0$ is as in Corollary \ref{cor:bounds} and $C$ is a constant independent of $\alpha$, $h$, $\eta$ and $\varepsilon$.
\end{theorem}
\begin{proof}
In the proof, we will make use of Young's inequality with weights, th\new{at} is\new{,}
\begin{equation}\label{eq:young}
xy \leq w x^q + \frac{1}{q'}\left(\frac{1}{qw}\right)^{\frac{q'}{q}} y^{q'},
\end{equation}
which holds true for all $x,y \geq 0$, $w > 0$, and $q \new{>} 1$ with $\frac1q + \frac1{q'} = 1$, see \cite[Eq. (3.1)]{hw22}.
	
Let us introduce the following notations for the occurring residuals:
\begin{align*}
	\res_{\infty} &:= \norm{F\left(\hat u_\alpha\right) - \gdag}{\L^\infty \left(\manifold\right)}, &\res_{1} &:= \norm{F\left(\hat u_\alpha\right) - \gdag}{\L^1 \left(\manifold\right)}, \\
	\res_{2} &:= \norm{F\left(\hat u_\alpha\right) - \gdag}{\L^2 \left(\manifold\right)},&\res_{\mathrm{d}} &:= \norm{F\left(\hat u_\alpha\right) - \gdag}{\Y_{\mathrm{d}}}.	
\end{align*}
In what follows, $C > 0$ denotes a generic constant independent of all asymptotic parameters $\eta, \varepsilon, h, \alpha$, the value of which might change from line to line. Furthermore, for simplicity, we ignore the dependency on $\rho$ in what follows, even though this could be made explicit.

Example \ref{ex:vsc} implies under Assumptions (A2) and (A3) the variational source condition \eqref{eq:vsc_L2}. Similar to the proof of \cite[Thm. 3.4]{hw22}, we obtain consequently
\begin{align*}
\frac14 \norm{\hat u_\alpha - \udag}{\L^2 \left(\Omega\right)}^2 \leq& \frac1\alpha\left(\norm{F_h\left(\udag\right) - \gobs}{\L^1 \left(\manifold_h\right)} - \norm{F_h \left(\hat u_\alpha\right) - \gobs}{\L^1 \left(\manifold_h\right)} \right) \\&+ C \left(\res_{2}^2\right)^{\frac{s}{s+a}}\\
\leq& \frac{1}{\alpha} \err \left(\hat u_\alpha\right) - \frac1\alpha \res_1 + C \res_2^{\frac{2s}{s+a}}\\
\leq& C \left[\frac{\varepsilon}{\alpha} + \frac{\eta}{\alpha} \res_\infty+ \frac{h^\omega}{\alpha}\res_{\mathrm{d}} - \frac1\alpha \res_1 +\res_2^{\frac{2s}{s+a}}\right],
\end{align*}
where we inserted the upper bound for $\err$ under the impulsive noise model from Theorem \ref{thm:admissibility} using Assumptions (A1) and (A5). Due to the $\left(\L^\infty\left(\manifold\right), C \cdot^{\frac{k}{d} - \frac1p}, \norm{\cdot - \cdot}{\L^2 \left(\Omega\right)}^2\right)$-smoothingness of $F$ due to Lemma \ref{lem:smoothingness} under Assumption (A4) and the assumed $\left(\Y_{\mathrm{d}}, \gamma, \norm{\cdot - \cdot}{\L^2 \left(\Omega\right)}^2\right)$-smoothingness of $F$ in Assumption (A6), the terms $\res_\infty$ and $\res_{\mathrm{d}}$ can be absorbed in $- \frac2{3\alpha} \res_1$ with additional error contributions $C\frac{\eta^{2 \vartheta+2}}{\alpha^2}$ and $C\frac{h^{2 \omega} \left(\gamma \left(3h^\omega\right)\right)^2}{\alpha^2}$, cf. the proof of Theorem 3.4 in \cite{hw22}. So we obtain
\begin{align*}
	\norm{\hat u_\alpha - \udag}{\L^2 \left(\Omega\right)}^2 
	\leq& C \left[\frac{\varepsilon}{\alpha} + \frac{\eta^{2 \vartheta+2}}{\alpha^2}+\frac{h^{2 \omega} \left(\gamma \left(3 h^\omega\right)\right)^2}{\alpha^2} - \frac1{3\alpha} \res_1 +\res_2^{\frac{2s}{s+a}}\right],
\end{align*}
Inspired by \cite[Prop. 2.2]{kwh16}, we now use
\begin{align*}
\res_2^2 \leq \res_\infty \res_1 &\leq \left(C \delta^{\frac{k}{d} - \frac1p} \norm{\hat u_\alpha - \udag}{\L^2 \left(\Omega\right)} + \frac1\delta \res_1 \right) \res_1\\
& = C \delta^{\vartheta} \norm{\hat u_\alpha - \udag}{\L^2 \left(\Omega\right)} \res_{1} + \frac1\delta \res_1^2,
\end{align*}
where we again used the $\left(\L^\infty\left(\manifold\right), C \cdot^{\frac{k}{d} - \frac1p}, \norm{\cdot - \cdot}{\L^2 \left(\Omega\right)}^2\right)$-smoothingness of $F$. Inserting this into our previous estimate, we obtain from $\left(x+y\right)^{\frac{s}{s+a}} \leq x^{\frac{s}{s+a}}+y^{\frac{s}{s+a}}$, that
\begin{align*}
\norm{\hat u_\alpha - \udag}{\L^2 \left(\Omega\right)}^2 \leq &C \left[\frac{\varepsilon}{\alpha} + \frac{\eta^{2 \vartheta+2}}{\alpha^2}+\frac{h^{2 \omega} \left(\gamma \left(3 h^\omega\right)\right)^2}{\alpha^2} - \frac1{3\alpha} \res_1 \right.\\
& \left. + \left(\delta^{\vartheta} \norm{\hat u_\alpha - \udag}{\L^2 \left(\Omega\right)} \res_{1}\right)^{\frac{s}{s+a}} + \left(\frac1\delta\right)^{\frac{s}{s+a}} \res_1^{\frac{2s}{s+a}}\right],
\end{align*}
The mixed term $\delta^{\vartheta} \norm{\hat u_\alpha - \udag}{\L^2 \left(\Omega\right)} \res_{1}$ can now be handled by Young's inequality \eqref{eq:young} with $w = \frac1\delta$ and $q = q' = 2$, which yields
\[
\left(\delta^{\vartheta} \norm{\hat u_\alpha - \udag}{\L^2 \left(\Omega\right)}\right) \res_{1} \leq \frac14\delta^{2\vartheta+1}\norm{\hat u_\alpha - \udag}{\L^2 \left(\Omega\right)}^2 +\frac{1}{\delta} \res_1^2.
\]
Inserting this once more into the previous estimate and using $\left(x+y\right)^{\frac{s}{s+a}} \leq x^{\frac{s}{s+a}}+y^{\frac{s}{s+a}}$ again, we find
\begin{align*}
	\norm{\hat u_\alpha - \udag}{\L^2 \left(\Omega\right)}^2 \leq &C \left[\frac{\varepsilon}{\alpha} + \frac{\eta^{2 \vartheta+2}}{\alpha^2}+\frac{h^{2 \omega} \left(\gamma \left(3 h^\omega\right)\right)^2}{\alpha^2} - \frac1{3\alpha} \res_1 \right.\\
	& \left. + \delta^{\frac{(2\vartheta+1)s}{s+a}}\norm{\hat u_\alpha - \udag}{\L^2 \left(\Omega\right)}^{\frac{2s}{s+a}}+2 \left(\frac1\delta\right)^{\frac{s}{s+a}} \res_1^{\frac{2s}{s+a}}\right],
\end{align*}
The terms $- \frac1{3\alpha} \res_1$ and $2 \left(\frac1\delta\right)^{\frac{s}{s+a}} \res_1^{\frac{2s}{s+a}}$ can now be handled by estimating
\[
-\frac1{3\alpha} \res_1 + 2\delta^{-\frac{s}{s+a}} \res_1^{\frac{2s}{s+a}} \leq \sup_{t \geq 0} \left[2\delta^{-\frac{s}{s+a}} t^{\frac{2s}{s+a}} - \frac{1}{3 \alpha} t\right] = C \left(\frac{\alpha^2}{\delta}\right)^{\frac{s}{a-s}}.
\]
Inserting this into the previous estimate yields
\begin{align*}
	\norm{\hat u_\alpha - \udag}{\L^2 \left(\Omega\right)}^2 \leq &C \left[\frac{\varepsilon}{\alpha} + \frac{\eta^{2 \vartheta+2}}{\alpha^2}+\frac{h^{2 \omega} \left(\gamma \left(3 h^\omega\right)\right)^2}{\alpha^2} \right.\\
	& \left. + \delta^{\frac{(2\vartheta+1)s}{s+a}}\norm{\hat u_\alpha - \udag}{\L^2 \left(\Omega\right)}^{\frac{2s}{s+a}}+\left(\frac{\alpha^2}{\delta}\right)^{\frac{s}{a-s}}\right],
\end{align*}
and handling the mixed term by Young \eqref{eq:young} with $w = \frac12$, $q = \frac{s+a}{s}$, $q' = \frac{s+a}{a}$, i.e.
\[
\delta^{\frac{(2\vartheta+1)s}{s+a}}\norm{\hat u_\alpha - \udag}{\L^2 \left(\Omega\right)}^{\frac{2s}{s+a}} \leq \frac12 \norm{\hat u_\alpha - \udag}{\L^2 \left(\Omega\right)}^2 + C \delta^{\frac{(2\vartheta+1)s}{a}},
\]
this implies by rearranging that
\begin{align*}
	\norm{\hat u_\alpha - \udag}{\L^2 \left(\Omega\right)}^2 \leq &C \left[\frac{\varepsilon}{\alpha} + \frac{\eta^{2 \vartheta+2}}{\alpha^2}+\frac{h^{2 \omega} \left(\gamma \left(3 h^\omega\right)\right)^2}{\alpha^2} + \delta^{\frac{(2\vartheta+1)s}{a}}+\left(\frac{\alpha^2}{\delta}\right)^{\frac{s}{a-s}}\right]
\end{align*}
for any $\delta \in (0,\delta_0)$. Optimizing this bound \new{with respect to} $\delta$ (at least in the limit $\alpha \searrow 0$ where the corresponding $\delta$ will satisfy $\delta < \delta_0$) yields \eqref{eq:bound_example}.
\end{proof}

\begin{remark}
	\begin{enumerate}
		\item[(a)] If $\Y_{\mathrm{d}} = \L^\infty \left(\manifold\right)$ in \eqref{eqs:discretization}, then Assumption (A6) of the above theorem is automatically satisfied under the other Assumptions according to Lemma \ref{lem:smoothingness}. 
		\item[(b)] If $F$ is at least $a$-smoothing for some $a \in \N$, $a > \frac{d}{2}$, i.e. \eqref{eq:a_smoothing_most} is complemented by \eqref{eq:a_smoothing_least}, then Assumption (A4) is automatically satisfied as well with $k = a$, $p = 2$. \new{As mentioned before, since $a > \frac{d}{2}$, the image of $F$ is contained in $\L^\infty \left(\manifold\right)$ and hence by boundedness of $\manifold$, the operator $F$ can be considered as a map into $\L^1\left(\manifold\right)$.}
		\item[(c)] Compared to the approach from \cite{kwh16}, where the $\res_\infty$-term in $\res_2^2 \leq \res_\infty \res_1$ is bounded using $\res_\infty \leq C$ uniformly, we get a different exponent of $\alpha$ in \eqref{eq:bound_example}. By the uniform bound one would obtain $\exp_{\mathrm{unif}}:=\frac{s}{a}$, whereas above we derived  $\exp_{\mathrm{fine}}:=\frac{2s(2\theta+1)}{(2\vartheta+1)(a-s)+a}$. In the situation of (b), i.e. if $F$ is at least $a$-smoothing in the sense of \eqref{eq:a_smoothing_least} and at most $a$-smoothing in the sense of \eqref{eq:a_smoothing_most}, we can choose $p =2$ and $k =a$, which results in $2\vartheta+1 = \frac{2a}{d}$, so that $\exp_{\mathrm{fine}} = \frac{4s}{2a-2s+d} > \exp_{\mathrm{unif}}$, \new{and in view of the additional requirement $a > d/2$ (see (b)),} our new bound is an improvement.
	\end{enumerate}	
\end{remark}

\section{Application to discretized noises and optimality in the Gaussian case}

\subsection{General setup}

Let us now apply the result from Theorem \ref{thm:example} to the situation of discretized noise. Therefore suppose that we have finitely many measurements $\gobs \in \R^n$ at hand, which are to be identified as
\begin{equation}\label{eq:discrete_data}
\gobs_i = \gdag \left(x_i\right) + \xi_i, \qquad i=1,\ldots,n,
\end{equation}
where $\xi_i$\new{,} $i =1,\ldots,n$\new{,} are i.i.d. centered random variables and $x_i \in \manifold$ are sampling points for $i=1,\ldots,n$. Here we assume w.l.o.g. that $\gdag$ can be interpreted as a continuous function, which is e.g. the case under Assumption (A4) from Theorem \ref{thm:example}. To include this setup into our previous model, we set
\[
\manifold_h := \left\{x_1,...,x_n\right\}, \qquad \Sigma = 2^{\manifold_h}, \qquad \nu \left(A\right):= \frac{\#A}{n} \left|\manifold\right|
\]
and
\[
F_h : D(F) \subset \X\to \L^1 \left(\manifold_h\right), \qquad F_h \left(u\right):= \left(\left(F(u)\right)\left(x_1\right),\ldots, \left(F(u)\right)\left(x_n\right)\right)^\top.
\]
Note that \eqref{eq:discretization2} is trivially fulfilled in this setup, and \eqref{eq:discretization1} holds true e.g. with $\Y_{\mathrm{d}} = \L^\infty\left(\manifold\right)$ as soon as $R(F)$ consists of sufficiently smooth functions, as is typically the case in inverse problems.

\new{Consider} a threshold $\lambda > 0$. If we define
\begin{align*}
	\manifold_{\mathrm{c}} = \left\lbrace i \in \left\lbrace 1, \ldots, n \right\rbrace  ~\big|~ \left|\xi_i \right| > \lambda\right\rbrace, 
\end{align*}
then this discrete data model fits into the setup of \eqref{eq:noise_model} with an impulsive noise function $\xi \in \L^1 \left(\manifold_h\right)$ with parameters
\begin{equation}\label{eq:epsilon_eta}
\varepsilon = \norm{\xi}{\L^1\left(\manifold_h\setminus \manifold_{\mathrm{c}}\right)} = \frac{\new{|\manifold|}}{n}\sum_{k\notin \manifold_{\mathrm{c}}} |\xi_k|, \qquad \eta = \nu\left(\manifold_{\mathrm{c}}\right)= \frac{\# \manifold_{\mathrm{c}} |\manifold|}{n}.
\end{equation}
Note that $\eta$ and $\varepsilon$ are random variables as well.

The moments of $\varepsilon$ and $\eta$ can now be bounded in terms of the cumulative distribution function $F_\xi$ of the random variables $\xi_i$, i.e.
\[
F_\xi \left(x\right) = \mathbb P \left[\xi_i \leq x\right], \qquad x \in \R,
\]
which is by the i.i.d. assumption independent of $i \in \{1,\ldots,n\}$:
\begin{lemma}\label{lem:moments}
	Let $\xi_i, i=1,\ldots,n$\new{,} i.i.d. centered random variables with cdf $F_\xi$. For the quantities $\varepsilon$ and $\eta$ defined in \eqref{eq:epsilon_eta}, it holds that
	\[
	\mathbb{E}[\varepsilon] \leq \new{|\manifold|}\sqrt{F_\xi \left(\lambda\right)-F_\xi\left(-\lambda\right)} \sqrt{\mathbb V \left[\xi_1\right]}\quad\text{and}\quad \mathbb{E}[\eta^{r}] \leq |\manifold|^{r} \left(1-F_\xi \left(\lambda\right)+F_\xi\left(-\lambda\right)\right)
	\]
	for any $r \geq 1$.
\end{lemma}
\begin{proof}
Recall $\manifold_{\mathrm{i}} := \manifold_h \setminus \manifold_{\mathrm{c}}$. We compute\new{,} using Cauchy-Schwarz\new{,} that
	\begin{multline*}
		\mathbb{E}[\varepsilon] = \E \left[\frac{\new{|\manifold|}}{n}\sum_{k\in \manifold_{\mathrm{i}}} |\xi_k|\right] = \frac{\new{|\manifold|}}{n} \E \left[\sum_{i=1}^n 1|_{\{|\xi_i| \leq \lambda\}} \left\|\xi_i\right\|\right]
		 = \frac{\new{|\manifold|}}{n} \sum_{i=1}^n\E \left[ 1|_{\{|\xi_i| \leq \lambda\}} \left|\xi_i\right|\right] \\= \new{|\manifold|}\E \left[1|_{\{|\xi_1| \leq \lambda\}} \left|\xi_1\right|\right]
		 \leq\new{|\manifold|} \sqrt{\E \left[ 1|_{\{|\xi_1| \leq \lambda\}}\right]} \sqrt{\E \left[\left|\xi_1\right|^2\right]}
	 = \new{|\manifold|}\sqrt{\mathbb P \left[\left|\xi_1\right| \leq \lambda\right]}\sqrt{\mathbb V \left[\xi_1\right]}\\
	= \new{|\manifold|}\sqrt{F_\xi \left(\lambda\right)-F_\xi\left(-\lambda\right)} \sqrt{\mathbb V \left[\xi_1\right]}
	\end{multline*}
	and using Hölder's inequality in the form
	\begin{multline*}
	\left(\sum_{i=1}^n y_i\right)^r = \left(\sum_{i=1}^n 1\cdot y_i\right)^r  \leq \left(\left(\sum_{i=1}^n 1\right)^{\frac{r-1}{r}} \left(\sum_{i=1}^n y_i^r\right)^{\frac1r} \right)^r\\ = \left(\sum_{i=1}^n 1\right)^{r-1} \left(\sum_{i=1}^n y_i^r\right) = n^{r-1} \sum_{i=1}^n y_i^r,
	\end{multline*}
	we find
	\begin{multline*}
		\mathbb{E}\left[\eta^{r}\right]= \E \left[\left(\frac{|\manifold|}{n}\#\manifold_{\mathrm{c}}\right)^{r}\right] 
		= \left(\frac{|\manifold|}{n}\right)^{r} \E \left[\left(\#\manifold_{\mathrm{c}}\right)^{r}\right] \\
		= \left(\frac{|\manifold|}{n}\right)^{r} \E \left[\left(\sum_{k=1}^{n} 1|_{\{|\xi_k| > \lambda\}}\right)^{r}\right] 
		\leq \left(\frac{|\manifold|}{n}\right)^{r} \E \left[n^{r-1}\sum_{k=1}^{n} 1|_{\{|\xi_k| > \lambda\}}^{r}\right] \\
		=  \frac{|\manifold|^{r}}{n} \E\left[\sum_{k=1}^n 1|_{\{|\xi_k| > \lambda\}}\right] 
		=  \frac{|\manifold|^{r}}{n} \sum_{k=1}^n \E\left[ 1|_{\{|\xi_k| > \lambda\}}\right] \\
		= |\manifold|^{r} \E\left[ 1|_{\{|\xi_1| > \lambda\}}\right]
		= |\manifold|^{r} \mathbb{P}\left(|\xi_1| > \lambda\right)\\
		= |\manifold|^{r} \left(1-F_\xi \left(\lambda\right)+F_\xi\left(-\lambda\right)\right).
	\end{multline*}
\end{proof}

This way, we can obtain error bounds for a variety of different discrete noise models.

\begin{remark}
Note that not only the moments of $\varepsilon$ and $\eta$ can be bounded, but it is also possible to derive concentration results by standard techniques. Therefore, all the following considerations would also allow for concentration results.
\end{remark}

\subsection{Gaussian case and optimal rates of convergence}

Let $\xi_i \stackrel{\mathrm{i.i.d.}}{\sim} \mathcal{N}\left(0,\sigma^2\right)$ for $i=1,\ldots,n$, i.e. consider the case of white Gaussian noise with variance $\sigma^2 > 0$. Let $\Phi: \R \to \R$ denote the standard cdf of $\mathcal N \left(0,1\right)$, i.e.
\[
\Phi(x) = \frac{1}{\sqrt{2\pi}} \int_{-\infty}^x \exp\left(-\frac{t^2}{2}\right) \,\mathrm d t, \qquad x \in \R,
\]
and recall that $\Phi \left(-x\right) = 1- \Phi(x)$ for any $x \in \R$.

For $\xi \sim \mathcal N\left(0,\sigma^2\right)$, the cumulative distribution function $F_\xi$ can be described in terms of $\Phi$: For such $\xi$ it holds $\frac{\xi}{\sigma} \sim \mathcal N \left(0,1\right)$, so that $F_{\frac{\xi}{\sigma}} = \Phi$, which implies
\[
F_\xi \left(x\right) = \mathbb P \left[\xi \leq x\right] = \mathbb P \left[\frac{\xi}{\sigma} \leq \frac{x}{\sigma}\right] = F_{\frac{\xi}{\sigma}}\left(\frac{x}{\sigma}\right) = \Phi \left( \frac{x}{\sigma}\right).
\]
Consequently, the result from Lemma \ref{lem:moments} reads as follows:
\begin{corollary}
Let $\xi_i \stackrel{\mathrm{i.i.d.}}{\sim} \mathcal N\left(0,\sigma^2\right)$ for $ i=1,\ldots,n$ with $\sigma\new{^2}> 0$. Then, for any $r \geq 1$, the moments of $\varepsilon$ and $\eta$ as defined in \eqref{eq:epsilon_eta} can be bounded by
\[
\mathbb{E}[\varepsilon] \leq \new{|\manifold|}\sigma\sqrt{1-2 \Phi \left(-\frac{\lambda}{\sigma}\right)} , \qquad \mathbb{E}[\eta^{r}] \leq 2|\manifold|^{r} \Phi \left(-\frac{\lambda}{\sigma}\right)
\]
for any $r \geq 1$.
\end{corollary}
\begin{proof}
	This follows directly from Lemma \ref{lem:moments} and
	\[
	F_\xi \left(\lambda\right)-F_\xi\left(-\lambda\right) = \Phi \left( \frac{\lambda}{\sigma}\right) - \Phi \left(-\frac{\lambda}{\sigma}\right) = 1-2 \Phi \left(-\frac{\lambda}{\sigma}\right) 
	\]
	or
	\[
	1-F_\xi \left(\lambda\right)+F_\xi\left(-\lambda\right) = 1-\Phi \left( \frac{\lambda}{\sigma}\right) + \Phi \left(-\frac{\lambda}{\sigma}\right) = 2 \Phi \left(-\frac{\lambda}{\sigma}\right),
	\]
	respectively.
\end{proof}
Inserting this into the result of Theorem \ref{thm:example} yields the following bounds:
\begin{corollary}
In the setting of Theorem \ref{thm:example}, grant Assumptions (A2)--(A6) and let $\gobs$ be given as in \eqref{eq:discrete_data} with $\xi_i \stackrel{\mathrm{i.i.d.}}{\sim} \mathcal N\left(0,\sigma^2\right)$ for $ i=1,\ldots,n$ with $\sigma> 0$. Then, for sufficiently small $\alpha, \sigma$ and $h$, any minimizer $\hat u_\alpha$ of \eqref{eq:tik} (if \new{it exists}) obeys the bound
\begin{align*}
\mathbb E \left[\norm{\hat u_\alpha -\udag}{\L^2 \left(\Omega\right)}^2\right] \leq& C \left[\frac{\sigma}{\alpha}\sqrt{1-2 \Phi \left(-\frac{\lambda}{\sigma}\right)}  + \frac{1}{\alpha^2}\Phi \left(-\frac{\lambda}{\sigma}\right)\right.\\& \left.+\frac{h^{2\omega}\left(\gamma \left(3 h^\omega\right)\right)^2}{\alpha^2}+ \alpha^{\frac{2s(2\vartheta+1)}{(2\vartheta+1)(a-s)+a}}\right],
\end{align*}
where $\vartheta = \frac{k}{d} - \frac1p > 0$ is as in Corollary \ref{cor:bounds}, $\lambda > 0$ is arbitrary, and $C$ is a constant independent of $\alpha$, $h$, $\eta$ and $\varepsilon$.
\end{corollary}

Let us now derive rates of convergence from this result. Therefore, we will assume that the discretization error, encoded in the term $h^{2\omega}\left(\gamma \left(3h^\omega\right)\right)^2$ is negligible compared to the statistical noise in the observations. Note that this is the typical situation in Inverse Problems, where $F$ is smoothing and hence both the exponent $\omega$ in \eqref{eqs:discretization} will be \new{comparatively} large and the function $\gamma$ will \new{decay} fast. In this situation, we only have to treat the bound
\[
	\mathbb E \left[\norm{\hat u_\alpha -\udag}{\L^2 \left(\Omega\right)}^2\right] \leq C \left[\frac{\sigma}{\alpha}\sqrt{1-2 \Phi \left(-\frac{\lambda}{\sigma}\right)}  + \frac{1}{\alpha^2}\Phi \left(-\frac{\lambda}{\sigma}\right) + \alpha^{\theta}\right],
\]
where we abbreviated $\theta := \frac{2s(2\vartheta+1)}{(2\vartheta+1)(a-s)+a}$. Let us first investigate what happens if $\alpha = \alpha \left(\sigma, \lambda\right)$ is chosen a priori and optimal:
\begin{lemma}
In the setting of Theorem \ref{thm:example}, grant Assumptions (A2)--(A6) and let $\gobs$ be given as in \eqref{eq:discrete_data} with $\xi_i \stackrel{\mathrm{i.i.d.}}{\sim} \mathcal N\left(0,\sigma^2\right)$ for $ i=1,\ldots,n$ with $\sigma> 0$. Suppose that the discretization error $h^{2\omega}\left(\gamma \left(3 h^\omega\right)\right)^2$ is negligible. Then, for an optimally chosen $\alpha = \alpha \left(\sigma,\lambda\right)$, any minimizer $\hat u_\alpha$ of \eqref{eq:tik} (if \new{it exists}) satisfies
\begin{equation}\label{eq:err_bound_gauss}
\mathbb E \left[\norm{\hat u_\alpha -\udag}{\L^2 \left(\Omega\right)}^2\right] \leq C \max\left\{ \left(\sigma\sqrt{1-2 \Phi \left(-\frac{\lambda}{\sigma}\right)}\right)^{\frac{\theta}{\theta+1}}, \left(\Phi \left(-\frac{\lambda}{\sigma}\right)\right)^{\frac{\theta}{\theta+2}} \right\},
\end{equation}
where $\theta= \frac{2s(2\vartheta+1)}{(2\vartheta+1)(a-s)+a},\vartheta = \frac{k}{d} - \frac1p$, and $\lambda>0$ is arbitrary.
\end{lemma}
\begin{proof}
We follow the proof of Thm. 4.1 in \cite{hw14}. With $A = \sigma\sqrt{1-2 \Phi \left(-\frac{\lambda}{\sigma}\right)}$ and $B = \Phi \left(-\frac{\lambda}{\sigma}\right)$ and any choice
\[
\alpha \sim \max \left\{A^{\frac1{\theta+1}}, B^{\frac1{\theta+2}}\right\},
\]
we obtain
\begin{multline*}
\mathbb E \left[\norm{\hat u_\alpha -\udag}{\L^2 \left(\Omega\right)}^2\right]  \leq C \left[\frac{A}{\alpha} + \frac{B}{\alpha^2} + \alpha^{\theta}\right] \leq C \left[\frac{A}{A^{\frac1{\theta+1}}} + \frac{B}{B^{\frac2{\theta+2}}} + \alpha^{\theta}\right] \\= C \alpha^\theta  \leq C\max \left\{A^{\frac{\theta}{\theta+1}}, B^{\frac{\theta}{\theta+2}}\right\},
\end{multline*}
which proves the claim.
\end{proof}

Note that this shows \new{(in principle)} a bound
\[
\mathbb E \left[\norm{\hat u_\alpha -\udag}{\L^2 \left(\Omega\right)}^2\right] \leq C \inf_{\lambda > 0}\max\left\{ \left(\sigma\sqrt{1-2 \Phi \left(-\frac{\lambda}{\sigma}\right)}\right)^{\frac{\theta}{\theta+1}}, \left(\Phi \left(-\frac{\lambda}{\sigma}\right)\right)^{\frac{\theta}{\theta+2}} \right\},
\]
which might however be hard to realize in practice, since the chosen regularization parameter $\alpha > 0$ would depend on $\lambda > 0$. Nevertheless, by choosing an appropriate $\lambda>0$, we find that the following convergence rate can be attained:
\begin{theorem} \label{thm:err_bound_gaussian}
In the setting of Theorem \ref{thm:example}, grant Assumptions (A2)--(A6) and let $\gobs$ be given as in \eqref{eq:discrete_data} with $\xi_i \stackrel{\mathrm{i.i.d.}}{\sim} \mathcal N\left(0,\sigma^2\right)$ for $ i=1,\ldots,n$ with $\sigma> 0$. Suppose that the discretization error $h^{2\omega}\left(\gamma \left(3 h^\omega\right)\right)^2$ is negligible. Then, for an optimally chosen $\alpha$, any minimizer $\hat u_\alpha$ of \eqref{eq:tik} (if \new{it exists}) satisfies
\[
\mathbb E \left[\norm{\hat u_\alpha -\udag}{\L^2 \left(\Omega\right)}^2\right]  = \mathcal O \left(\sigma^{\frac{\theta}{\theta+1}}\right) \qquad\text{as}\qquad\sigma \searrow 0,
\]
where $\theta= \frac{2s(2\vartheta+1)}{(2\vartheta+1)(a-s)+a},\vartheta = \frac{k}{d} - \frac1p$.
\end{theorem}
\begin{proof}
We aim to choose $\lambda > 0$ such that both terms in the maximum in \eqref{eq:err_bound_gauss} equal each other. Since we want the right-hand side to go to $0$ as $\sigma \searrow 0$, this requires $z = \Phi \left(-\frac{\lambda}{\sigma}\right) \to 0$ as well. Consequently, the $\sqrt{1-2z}$-term in the first argument of the maximum can just be ignored, as it will be bounded for $z \to 0$. Since, for every $\sigma > 0$, there exists a $\lambda > 0$ (with the asymptotics $\lambda\to \infty$ as $\sigma \searrow 0$) with $\sigma^{\frac{\theta}{\theta+1}} = z = \Phi \left(-\frac{\lambda}{\sigma}\right)$, this shows that both terms in the maximum \eqref{eq:err_bound_gauss} can be equa\new{ted} (asymptotically), yielding the claimed bound.
\end{proof}

Let us now specify the result from Theorem \ref{thm:err_bound_gaussian} in case of an $a$-smoothing operator, i.e. where \eqref{eq:a_smoothing_most} and \eqref{eq:a_smoothing_least} hold true simultaneously. The convergence rate derived in the Corollary below coincides with \eqref{eq:rate_opt}, i.e. is order-optimal in this situation, which proves that the usage of $\L^1$ data fitting does not deteriorate the possible rate of convergence.

\begin{corollary}\label{cor:gauss_opt}
In the setting of Theorem \ref{thm:example}, grant Assumptions (A3), (A5) and (A6) and assume that $F$ is $a$-times smoothing for some $a \in \N$, $a > d/2$ as in \eqref{eq:a_smoothing_most} and \eqref{eq:a_smoothing_least}. Let $\gobs$ be given as in \eqref{eq:discrete_data} with $\xi_i \stackrel{\mathrm{i.i.d.}}{\sim} \mathcal N\left(0,\sigma^2\right)$ for $ i=1,\ldots,n$ with $\sigma> 0$. Suppose that the discretization error $h^{2\omega}\left(\gamma \left(3 h^\omega\right)\right)^2$ is negligible. Then, for an optimally chosen $\alpha$, any minimizer $\hat u_\alpha$ of \eqref{eq:tik} (if \new{it exists}) satisfies
\[
\left(\mathbb E \left[\norm{\hat u_\alpha -\udag}{\L^2 \left(\Omega\right)}^2\right]\right)^{\frac12} = \mathcal O \left(\sigma^{\frac{2s}{2a+2s+d}}\right).
\]
\end{corollary}
\begin{proof}
As mentioned before, (A4) with $k = a, p=2$ is a direct consequence of \eqref{eq:a_smoothing_least}. With these parameters, we obtain $\vartheta = \frac{a}{d} - \frac12$, so that $2\vartheta+1 = \frac{2a}{d}$ and $\theta = \frac{4s}{2a-2s+d}$. This yields	
\[
\frac{\theta}{\theta+1} = \frac{4s}{2a+2s+d},
\]
which proves the claim.
\end{proof}
%
%\subsection{Laplace case}
%
%Let $\xi_i \stackrel{\mathrm{i.i.d.}}{\sim} \text{Lap}\left(0,\sigma\right)$ for $i=1,\ldots,n$, i.e. consider the case of centered Laplace distributed variables with scaling parameter $\sigma > 0$. Note, that Laplacian noise is of special interest, since the negative log-likelihood in finite dimensions is just given by a $1$-norm.
%
%For $\xi\sim\text{Lap}\left(0,\sigma\right)$, the cdf $F_\xi$ is given by
%\[
%F_\xi\left(x\right) = \frac12 + \frac12 \sgn(x) \left(1-\exp\left(-\frac{|x|}{\sigma}\right)\right),
%\]
%where $\sgn$ denotes the usual sign-function, i.e. $\sgn(x) = 1$ if $x > 0$, $\sgn(x) = -1$ if $x < 0$ and $\sgn(0) = 0$. Consequently, the result from Lemma \ref{lem:moments} reads as follows:
%
%\todo[inline]{Derive bounds exactly as in Section 4.1}

\section{Numerical simulations}

In this \new{s}ection, we will now present an Algorithm for the solution of \eqref{eq:tik} in case of a bounded linear operator $T :\L^2 \left(\Omega\right) \to \L^{\new{1}} \left(\manifold\right)$ with $\mathcal R(u) = \norm{u}{\L^2 \left(\Omega\right)}^2$, i.e.
\begin{align}\label{eq:primal}
    \widehat{u}_\alpha \in \argmin_{u \in \L^1(\manifold)} \left[ \norm{Tu- \gobs}{\L^1(\manifold)} + \alpha\norm{u}{\L^2(\Omega)}^2 \right].
\end{align}
Note that we have set $F_h = F = T$ here. In this linear setup with quadratic penalty, existence and uniqueness of a minimizer can readily be derived as in \cite{s08}.

After deriving the minimization algorithm, we will also present numerical simulations for a test problem to validate the theoretical results from Section 4.

\subsection{Numerical algorithm}

In \cite{hw14}, the problem \eqref{eq:primal} has been solved by means of the dual problem, corresponding to a quadratic minimization problem under box constraints, which was then solved using the quadprog algorithm. This method has several advantages, including computable bounds for the accuracy by means of the primal-dual gap, but relying on quadprog makes it slow and \new{is} especially not applicable for large-scale problems, where the matrix corresponding to $F_h$ cannot be stored or set up explicitly. 

We will therefore apply ADMM schemes to \eqref{eq:primal} to solve the primal problem directly, cf. \cite[Sect. 15.2]{b17}. Substituting $ v=Tu-\gobs$, the problem \eqref{eq:primal} is equivalent to
\begin{align*}
    \argmin_{u,v} \left[ \norm{v}{\L^1} + \alpha\norm{u}{\L^2}^2 \right] \quad  \text{s.t.} \quad Tu-v=\gobs.
\end{align*}
Here and in what follows, we omit the domains in the definition of $\L^1$ and $\L^2$ for brevity. Applying the standard ADMM to this problem, we obtain the update scheme
\begin{align}
u_{k+1} &= \argmin_u \left[ \alpha \norm{u}{\L^2}^2 + \frac{\rho}{2} \norm{Tu-v_k-\gobs+ \frac{1}{\rho} \mu_k}{\L^2}^2 \right], \tag{ADMM1}\label{eq:admm1}\\
v_{k+1} &= \argmin_v \left[ \norm{v}{\L^1} +  \frac{\rho}{2} \norm{T u_{k+1} -v-\gobs + \frac{1}{\rho} \mu_k}{\L^2}^2\right], \tag{ADMM2}\label{eq:admm2}\\
\mu_{k+1} &= \mu_k + \rho \left(T u_{k+1}- v_{k+1}-\gobs\right) \tag{ADMM3}\label{eq:admm3}
\end{align}
Note that \eqref{eq:admm1} can be solved explicitly by means of the necessary and sufficient first-order optimality condition, yielding
\begin{align}\tag{ADMM1'}
    u_{k+1} = \left(2 \alpha I + \rho T^\ast T\right) ^{-1} \left(\rho T^\ast (v_k+\gobs) - T^\ast \mu_k\right).
\end{align}
This problem can be solved matrix-free e.g. by using the CG algorithm. Furthermore, $v_{k+1}$ in \eqref{eq:admm2} can be determined by soft-thresholding. If we introduce
\[
\psi^{\mathrm{soft}}_\tau \left(x\right) := \begin{cases} x-\tau & \text{if } x \geq \tau, \\ 0 & \text{if } \abs{x} \leq \tau, \\ x+\tau & \text{if } x < - \tau \end{cases}
\]
for $x \in \R$ and a parameter $\tau > 0$, then \eqref{eq:admm2} corresponds to
\begin{align}\tag{ADMM2'}\label{eq:admm2rev}
    v_{k+1} =  \psi^{\mathrm{soft}}_{\frac1\rho} \left(\frac{\mu_k}{\rho}+ Tu_{k+1}-\gobs\right),
\end{align}
where the application is \new{performed} point-wise.

Even though (ADMM) can be implemented matrix-free, it will still be computationally demanding due to the inversion in (ADMM1'). Therefore, we consider another ADMM-based method, which avoids this issue. It arises from adding a quadratic proximity term to the objective in \eqref{eq:admm1}. Given a coercive operator $G$ on $\L^2 \left(\Omega\right)$, such that $\norm{u}{G}^2 := \left\langle G u, u\right\rangle_{\L^2}$ defines a scalar product on $\L^2 \left(\Omega\right)$, we replace \eqref{eq:admm1} by
\[
\tilde u_{k+1} = \argmin_u \left[\alpha \norm{u}{\L^2}^2 + \frac{\rho}{2} \norm{Tu-v_k-\gobs+ \frac{1}{\rho} \mu_k}{\L^2}^2 + \frac{1}{2} \norm{u - \tilde u_k}{G}^2 \right].
\]
If now $G$ is chosen as $\beta \text{id}_{\L^2} - \rho T^*T$, where $\beta > \rho\lambda_{\max}\left(T^*T\right)$, then straightforward computations show that
\begin{align*}
&\alpha \norm{u}{\L^2}^2 + \frac{\rho}{2} \norm{Tu-v_k-\gobs+ \frac{1}{\rho} \mu_k}{\L^2}^2 + \frac{1}{2} \norm{u - \tilde u_k}{G}^2\\
=&\alpha \norm{u}{\L^2}^2 + \frac{\rho}{2} \norm{Tu}{\L^2}^2 - \rho \left\langle Tu, v_k + \gobs - \frac1\rho \mu_k\right\rangle_{\L^2} + \frac12\left\langle G \left(u-\tilde u_k\right), u-\tilde u_k\right\rangle_{\L^2} + c, \\
=& \left(\alpha + \frac\beta2\right) \norm{u}{\L^2}^2 - \rho \left\langle Tu, v_k + \gobs - \frac{1}{\rho} \mu_k - T \tilde u_k\right\rangle_{\L^2} - \beta \left\langle u, \tilde u_k\right\rangle_{\L^2}+c\\
=& \left(\alpha + \frac\beta2\right) \norm{u}{\L^2}^2 - \rho \left\langle u, T^* \left(v_k + \gobs - \frac{1}{\rho} \mu_k - T \tilde u_k\right) + \beta \tilde u_k\right\rangle_{\L^2}+c,
\end{align*}
where $c$ denotes a constant independent of $u$ which might change from line to line. Therefore, again exploiting the necessary and sufficient first-order optimality condition, $\tilde u_{k+1}$ is given as
\begin{align*}
\tilde u_{k+1} &= \argmin_u \left[\left(\alpha + \frac\beta2\right) \norm{u}{\L^2}^2 - \rho \left\langle u, T^* \left(v_k + \gobs - \frac{1}{\rho} \mu_k - T \tilde u_k\right) + \beta \tilde u_k\right\rangle_{\L^2} \right]\\
& = \left(2 \alpha + \beta \right)^{-1} \rho \left(T^* \left(v_k + \gobs - \frac{1}{\rho} \mu_k - T \tilde u_k\right) + \beta \tilde u_k\right)\\
& = \left(1+ 2 \frac{\alpha}{\beta}\right)^{-1} \left( \tilde u_k - \frac\rho\beta T^* \left(T \tilde u_k - v_k - \gobs + \frac1\rho \mu_k\right)\right).
\end{align*}
Note that this can also be interpreted as linearizing the quadratic term in \eqref{eq:admm1} and adding a penalty, which is why this method is called the alternating direction linearized proximal method of multipliers (AD-LPMM). The resulting algorithm is as follows:
\begin{align}
	\tag{AD-LPMM1} \label{eq:adlpmm1}
	\tilde u_{k+1} &= \left(1+ 2 \frac{\alpha}{\beta}\right)^{-1} \left( \tilde u_k - \frac\rho\beta T^* \left(T \tilde u_k - v_k - \gobs + \frac1\rho \mu_k\right)\right),\\
	\tag{AD-LPMM2}\label{eq:adlpmm2}v_{k+1}& =  \psi^{\mathrm{soft}}_{\frac1\rho} \left(\frac{\mu_k}{\rho}+ T\tilde u_{k+1}-\gobs\right)\\
	\tag{AD-LPMM3}\label{eq:adlpmm3}\mu_{k+1} &= \mu_k + \rho \left(T \tilde u_{k+1} - v_{k+1}-\gobs\right).
	\end{align}
Compared to \eqref{eq:admm1}, \eqref{eq:adlpmm1} does not require the solution of a linear system and can hence be computed much faster. Besides this, the second and third step in ADMM and AD-LPMM, i.e. \eqref{eq:admm2rev} and \eqref{eq:adlpmm2} as well as \eqref{eq:admm3} and \eqref{eq:adlpmm3}, coincide. We emphasize that $\beta$ has to be chosen sufficiently large to ensure $\beta > \rho\lambda_{\max}\left(T^*T\right)$. In practice, this can be achieved by a few steps of the power iteration. For a convergence analysis of ADMM and AD-LPMM we refer to \cite{b17}. 
	
\subsection{Test problem}

For our numerical tests, we consider the problem from \cite{hw14}, i.e. the linear integral operator $T : \L^2 \left(\left[0,1\right]\right) \to \L^2 \left(\left[0,1\right]\right)$ given by
\[
\left(Tu\right)\left(x\right) = \int_0^1 k \left(x,y\right) u(y) \,\mathrm d y, \qquad u \in \L^2 \left(\left[0,1\right]\right)
\]
with kernel
\[
k\left(x,y\right) = \min \left\{ x \left(1-y\right), y \left(1-x\right)\right\}, \qquad x,y \in [0,1].
\]
Note that $\left(Tu\right)^{\prime\prime} = - u$ for any function $u \in \L^2 \left(\left[0,1\right]\right)$, so that estimating $u$ from $Tu$ corresponds to estimating the second derivative of $Tu$. Obviously, $T$ is $2$-smoothing, such that \eqref{eq:a_smoothing_most} and \eqref{eq:a_smoothing_least} are both satisfied with $a = 2$. \new{Since $a > d/2 = 1/2$, we can also consider $T$ as an operator mapping into $\L^1 \left([0,1]\right)$.}

To discretize $T$, we use the composite midpoint rule with $n$ equidistant evaluation points $x_1 = \frac1{2n}, x_2 = \frac3{2n}, \ldots,x_n =  \frac{2n-1}{2n}$, i.e.
\begin{equation}\label{eq:discretization}
\left(Tu\right)\left(x\right) = \int_0^1 k \left(x,y\right) u(y) \,\mathrm d y \approx  \frac 1n \sum_{i=1}^n k\left(x,x_i\right) u \left(x_i\right).
\end{equation}
%where the parameter $h = \frac1n$. This yields a discretization $T_{\mathrm{discrete}} : \L^2 \left(\left[0,1\right]\right) \to \L^2 \left(\left[0,1\right]\right)$, and by taking samples at the same grid points, we obtain the approximation operator $T_h: \L^2 \left(\left[0,1\right]\right) \to  \R^n$ in as
%\[
%T_h u := \left(\left(T_{\mathrm{discrete}} u\right) \left(x_1\right), \ldots,  \left(T_{\mathrm{discrete}} u\right) \left(x_n\right)\right)^\top  = \left(\frac 1n \sum_{i=1}^n k\left(x_j,x_i\right) u \left(x_i\right)\right)_{j=1}^n.
%\]
%It can be shown by standard results from numerical integration that \eqref{eqs:discretization} are satisfied and that the resulting error $h^\omega$ can be neglected compared to reasonable noise levels $\sigma$. Note, that in our discretized setup this error does not play a role at all since also $\hat u_\alpha$ can be computed only in a discretized fashion. This additional source of error would only pop up if we would compare different discretization levels.

As true solution, we choose
\[
\udag \left(x\right) = \begin{cases} x &\text{if } 0 \leq x \leq \frac12, \\ 1-x & \text{if } \frac12 \leq x \leq 1. \end{cases}
\]
Then $u \in \H^s \left(\left[0,1\right]\right)$ for any $s < \frac32$ as shown in \cite{lw20}, and thus all assumptions from Theorem \ref{thm:example} are satisfied with $F = F_h = T$. The resulting convergence rate in Theorem \ref{thm:err_bound_gaussian} is $\mathcal O \left(\sigma^{\frac{2s}{2s+4+1}}\right)$. In the limiting case $s = \frac32$, which cannot be handled by our analysis, but which can be reached by replacing Sobolev with Besov spaces, we obtain the optimal rate of convergence $\mathcal O \left(\sigma^{\frac{3}{8}}\right)$ for the function $\udag$, see \cite{hm19} and \cite[Rem. 2.3]{hw22}.

To avoid an inverse crime, the function $\gdag = T\udag$ is computed analytically and independent of the discretization.

\subsection{Results}

In the following, we depict both reconstructions and simulated rates of convergence for two settings: standard $\L^2$ data fitting, i.e. 
\begin{equation}\label{eq:tikl2}
	\hat u_\alpha= \argmin_{u \in\L^2} \left[ \norm{Tu - \gobs}{\L^2}^2 + \alpha \norm{u}{\L^2}^2\right],
\end{equation}
and \eqref{eq:primal} itself. Furthermore, we present all these simulations both carried out using (ADMM) and (AD-LPMM).

To simulate the rate of convergence, we use $n = 257$ in \eqref{eq:discretization} and $M = 100$ Monte Carlo runs. In each run, we simulate data $\gobs$ as in \eqref{eq:discrete_data} and perform \eqref{eq:primal} both with ADMM and AD-LPMM as well as \eqref{eq:tikl2}. Therein, we perform $K = 10\new{,}000$ iterations both in ADMM and AD-LPMM. The regularization parameter $\alpha$ is varied over $40$ different values, which are chosen logarithmically equispaced between $10^{-8}$ and $10^{-2}$, separately for each set of data. Afterwards, we determine that value of $\alpha$, which leads \new{to} the minimal error $\frac1M \sum_{m=1}^M \norm{\hat u_\alpha^i - \udag}{\L^2}^2$, where $i$ denotes the reconstruction in the $i$th run, $i =1,...,M$. This way, we simulate the optimal $\alpha$ as in Theorem \ref{thm:err_bound_gaussian}. The above procedure is repeated for different values of $\sigma$, and the corresponding optimal error is plotted against $\sigma$ in Figure \ref{fig:results} below.

\setlength{\fwidth}{11cm}
\setlength{\fheight}{5cm}

\begin{figure}[!htb]
	\centering
	\footnotesize
	\begin{tikzpicture}[baseline]
	\begin{axis}[%
		width=\fwidth,
		height=\fheight,
		scale only axis,
		xmode=log,
		ymode = log,
		xmin=1e-6,
		xmax=1e-1,
		xminorticks=true,
		ymin=1e-2,
		ymax=10,
		xlabel = $\sigma$,
		ylabel = $\sqrt{\frac1M \sum_{i=1}^M \norm{\hat u_{\alpha_{\mathrm{opt}}}^i - \udag}{\L^2}^2}$,
		legend pos=north west,
		legend cell align=left
		]
		\addplot[color=red] table [x index=0, y index=1] {adlpmm.dat};
		\addlegendentry{$\L^1$: \eqref{eq:tik} solved by AD-LPMM} 
		\addplot[color=green] table [x index=0, y index=1] {admm.dat};
		\addlegendentry{$\L^1$: \eqref{eq:tik} solved by ADMM} 
		\addplot[color=blue] table [x index=0, y index=1] {L2.dat}; 
		\addlegendentry{$\L^2$: \eqref{eq:tikl2}} 
		\addplot[color = black, domain=5e-5:5e-2, samples=200] {\x^(3/8)};
		\addlegendentry{$\mathcal O\left(\sigma^{\frac38}\right)$}
	\end{axis}
	\end{tikzpicture}%
	\caption{Simulations results for \eqref{eq:tikl2} and \eqref{eq:primal}, the latter solved either by ADMM or AD-LPMM. The different lines in color depict the error $\left(\mathbb E \left[\norm{\hat u_\alpha -\udag}{\L^2}^2\right]\right)^{\frac12}$ approximated by $M = 100$ Monte Carlo runs, where the optimal $\alpha_{\mathrm{opt}}$ is chosen out of a set of $40$ candidates $\alpha \in \left[10^{-8}, 10^{-2}\right]$. The black line illustrates the optimal rate of convergence $\mathcal O \left(\sigma^{\frac38}\right)$.}
	\label{fig:results}
\end{figure}
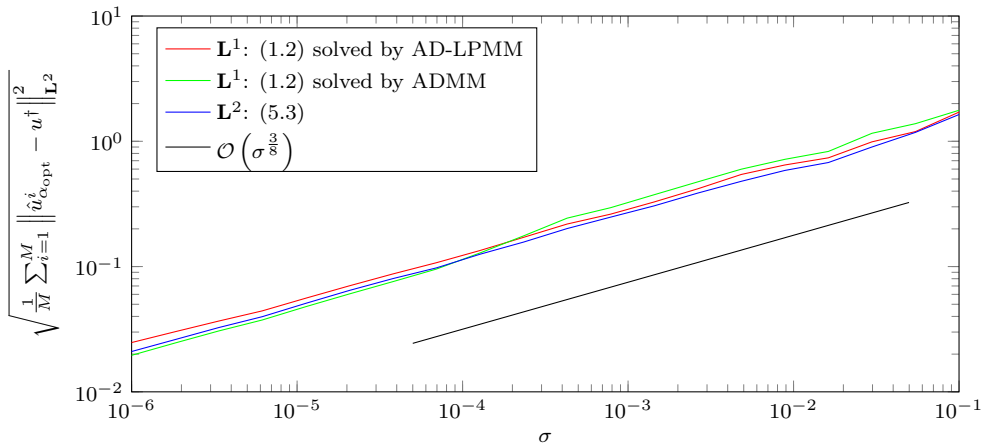

The results in Figure \ref{fig:results} support the theory from Theorem \ref{thm:err_bound_gaussian} and Corollary \ref{cor:gauss_opt}. There is visually no difference between the observed convergence rate for standard $\L^2$ Tikhonov regularization \eqref{eq:tikl2}, which is well-known to attain the optimal rate $\mathcal O \left(\sigma^{\frac38}\right)$ of convergence, and \eqref{eq:primal}, no matter if solved via ADMM and AD-LPMM. 

\section{Conclusion}

In this paper, we have investigated $\L^1$ data fitting for Inverse Problems with discretized noises. We have derived error bounds in a general setting, including noise in the forward operator, under verifiable conditions. Applied to the situation of $a$-smoothing operators and Gaussian white noise, this yields order-optimal rates of convergence. Therefore, the usage of $\L^1$ data fitting comes at no price, at least asymptotically. We have also discussed possible algorithms for the solution of the minimization problem, which allow \new{computing} the corresponding minimizers even in large-scale real world applications. In view of the increased robustness of $\L^1$ data fitting compared to standard $\L^2$ based approaches, it seems therefore reasonable to use $\L^1$ data fitting more often.

\new{A}n interesting question for future research, which was not treated in this paper, is how to choose the parameter $\alpha$ in \eqref{eq:tik}. As discussed in Section 4, the derived order-optimal convergence rates are only valid for an optimal choice of $\alpha$ depending on the thresholding parameter $\lambda$, which itself is chosen depending on unknown quantities. Therefore, an a posteriori parameter choice rule would be of immediate interest.

Future research directions could also include further development of the algorithms, including computable error bounds after a finite number of iterations and reliable stopping criteria. Simulations for further and especially larger problems are desirable.

\bibliographystyle{amsalpha}

\bibliography{l1_data_fitting}

\end{document}